\documentclass[12pt,a4paper]{article}

\usepackage{a4wide}

\usepackage[intlimits]{amsmath}
\usepackage{amssymb,amsthm}

\usepackage[cp1251]{inputenc}
\usepackage[english,ukrainian]{babel}
\begin{document}
\vskip 20mm {\large

\vskip 5mm \leftline{{\bf I.V.Samoilenko}} \vskip 5mm

\leftline{\normalsize Institute of Mathematics NASU, Tereschenkivska
str., 3, 01601, Kyiv, Ukraine} \vskip 5mm \leftline{\normalsize
isamoil@imath.kiev.ua}
\vskip 5mm

\leftline{\large \bf Weak convergence of Markovian random evolution}

\leftline {\large \bf in a multidimensional space.}

\vskip 15mm

{\normalsize }

\textbf{Abstract}. We study Markovian symmetry and non-symmetry
random evolutions in $\mathbf{R}^n$. Weak convergence of Markovian
symmetry random evolution to Wiener process and of Markovian
non-symmetry random evolution to a diffusion process with drift is
proved using problems of singular perturbation for the generators of
evolutions. Relative compactness in
$\mathbf{D}_{\mathbf{R}^n\times\Theta}[0,\infty)$ of the families of
Markovian random evolutions is also shown.



\newpage

\vskip 10mm

\section{Introduction}

\textit{Markov symmetry random evolutions} (MSRE) in spaces of
different dimensions were studied in the works of M.Kac \cite{Kac},
M.Pinsky \cite{Pin}, E.Orzingher (e.g., \cite{Or1, Or2}), A.F.Turbin
and A.D.Kolesnik (e.g., \cite{KoTu1, KoTu2})(see also \cite{Sam1}
for other references). Symmetry in this sense should be regarded as
uniform stationary distribution of switching at a symmetrical
structure in  $\mathbf{R}^n$, for instance at $n+1$-hedron
\cite{Sam1}, or at a unit sphere \cite{Kol2}.

Weak convergence of distributions of MSRE is also studied in some of
these works, namely convergence in $\mathbf{R}^2$ and $\mathbf{R}^3$
 was proved by A.D.Kolesnik in \cite{Kol1, Kol2}.

We should note that the problem of weak convergence of random walks
(partially, similar to MSRE) was studied by many authors. Among the
most interesting works we may point \cite{BoBo, Bo, FoKo, Fo}. Large
bibliography as for this branch could be found in \cite{BoBo}. The
methods, proposed in these works let us solve a wide range of
problems connected with convergence of random walks, but do not let
to obtain limit process to be averaged by the stationary measure of
switching process.

Such averaging may be found in the works of V.V.Anisimov and his
students (see \cite{An} and references therein), but here the
averaging by the stationary measure is one of the conditions,
proposed for the pre-limit process in the corresponding theorem.

In this work MSRE in $\mathbf{R}^n$ is studied using the methods,
proposed in \cite{korlim}. We find a solution of singular
perturbation problem for the generator of the evolution and thus the
averaging by a stationary measure of switching process is obtained
as a corollary of this solution. At the second stage we prove
relative compactness of the family of MSRE. This method let us show
weak convergence of the process of MSRE to the Wiener process in
$\mathbf{R}^n$.

The difference in the methods may be easily seen by the analysis of
the papers \cite{An1} and \cite{Sam3}, where similar problems are
studied.

In the sections 4 and 5 we use the method proposed to prove weak
convergence of \textit{Markov non-symmetry random evolution} (MNRE)
in $\mathbf{R}^n$. The distinction of this model is that the limit
process is a diffusion process with deterministic drift.

\section{Description of MSRE}

We study a particle in the space $\mathbf{R}^n$, that starts at
$t=0$ from the point $x=(x_i, i=\overline{1,n})$. Possible
directions of motion are given by the vectors
$$s(\theta)=(\cos\theta_1,
\sin\theta_1\cos\theta_2,\sin\theta_1\sin\theta_2\cos\theta_3,...,\sin\theta_1...\sin\theta_{n-2}\cos\theta_{n-1},$$
$$\sin\theta_1...\sin\theta_{n-2}\sin\theta_{n-1}), \theta_{n-1}\in
[0,2\pi), \theta_i\in [0,\pi), i=\overline{1,n-2}.$$ These vectors
have initial point in the center of the unit $n$-dimentional sphere
$S_n$ and the terminal point at its surface. Challenge of every next
direction is random and its time is distributed by Poisson. Thus,
the switching process is Poisson one with intensity
$\lambda=\varepsilon^{-2}$. The velocity of particle's motion is
fixed and equals $v=c\varepsilon^{-1}$, where $\varepsilon$ - is a
small parameter, $\varepsilon\to 0 \hskip 2mm (\varepsilon>0)$.

Let's define a set $\Theta=\{\theta:s(\theta)\in S_n\}.$ The
switching Poisson process is $\theta_t^{\varepsilon}\in \Theta.$

\textbf{Definition 1:} \textit{Markov symmetry random evolution}
(MSRE) is the process $\xi_t^{\varepsilon}\in \mathbf{R}^n$, that is
given by:
$$\xi_t^{\varepsilon}:=x+v\int_0^ts(\theta_{\tau}^{\varepsilon})d\tau.$$

Easy to see that when $\varepsilon\to 0 \hskip 2mm (\varepsilon>0)$
the velocity of the particle and intensity of switching decrease.
Our aim is to prove weak convergence of MSRE to the Wiener process
when $\varepsilon\to 0$. The main method is solution of singular
perturbation problem for the generator of MSRE. Let's describe this
generator.

Two-component Markov process $(\xi_t^{\varepsilon},
\theta_t^{\varepsilon})$ at the test-functions
$\varphi(x_1,...,x_n;\theta)\in C^{\infty}_0(\mathbf{R}^n\times
\Theta)$ may be described by a generator (see, e.g. \cite{Pin})
\begin{equation} \label{GenEp} L^{\varepsilon}\varphi(x_1,...,x_n;\theta):=\lambda Q\varphi(\cdot;\theta)+v S(\theta)\varphi(x_1,...,x_n;\cdot)=$$ $$\varepsilon^{-2}Q\varphi(\cdot;\theta)+\varepsilon^{-1}cS(\theta)\varphi(x_1,...,x_n;\cdot),\end{equation}
where
$$S(\theta)\varphi(x_1,...,x_n;\cdot):=-\left(s(\theta),\nabla\right)\varphi(x_1,...,x_n;\cdot),$$
here $\nabla\varphi=\left(\partial\varphi/\partial x_i, 1\leq i\leq
n\right)$,
$$Q\varphi(\cdot;\theta):=(\Pi-I)\varphi(\cdot;\theta):=\frac{1}{N}\int_{S_n}\varphi(\cdot;\theta)\mu(d\theta)-\varphi(\cdot;\theta),$$
here $N=(2\pi)^{n/2}\frac{1}{2\cdot 4\cdot...\cdot(n-2)}$ for even
$n$, and $N=(2\pi)^{(n-1)/2}\frac{2}{3\cdot 5\cdot...\cdot(n-2)}$
for odd $n$; $\mu(d\theta)$ - is the element of volume of the sphere
$S_n$, that is equal to
$$\mu(d\theta):=\sin^{n-2}\theta_1\sin^{n-3}\theta_2...\sin\theta_{n-2}d\theta_1...d\theta_{n-1}.$$

Using well-known formula
$$\int_0^{\pi}sin^{2m}d\theta=2\cdot\frac{1\cdot2\cdot...\cdot (2m-1)}{2\cdot4\cdot...\cdot2m}\cdot\frac{\pi}{2},
\int_0^{\pi}sin^{2m+1}d\theta=2\cdot\frac{2\cdot4\cdot...\cdot
2m}{1\cdot3\cdot...\cdot(2m+1)}$$ we may see
$$\int_{S_n}\mu(d\theta)=N.$$

Operator
$\Pi\varphi(\cdot;\theta):=\frac{1}{N}\int_{S_n}\varphi(\cdot;\theta)\mu(d\theta)$
is the projector at the null-space of reducibly-invertible operator
$Q$, because by definition it transfers functions to constants, but
constants to itself.

For $\Pi$ we have: $$Q\Pi=\Pi Q=0.$$

Potential operator $R_0$ may be defined by:
$$R_0:=\Pi-I.$$

This operator has the property: $$R_0Q=QR_0=I-\Pi,$$ thus it is
inverse for $Q$ in the range of $Q$, but for the function $\phi$
from null-space of $Q$ we have
$$R_0\phi=0.$$

Solution of singular perturbation problem in the series scheme with
the small series parameter $\varepsilon\to 0 (\varepsilon>0)$ (see
\cite{korlim}) for reducibly-invertible operator $Q$ and perturbed
operator $Q_1$ consists in the following.

We should find a vector
$\varphi^{\varepsilon}=\varphi+\varepsilon\varphi_1+\varepsilon^2\varphi_2$
and a vector $\psi$, that satisfy asymptotic representation
$$[\varepsilon^{-2}Q+\varepsilon^{-1}Q_1]\varphi^{\varepsilon}=\psi+\varepsilon\theta^{\varepsilon}$$
with the vector $\theta^{\varepsilon}$, that is uniformly bounded by
the norm and such that $$||\theta^{\varepsilon}||\leq C,
\varepsilon\to 0.$$

The left part of the equation may be rewritten
$$[\varepsilon^{-2}Q+\varepsilon^{-1}Q_1](\varphi+\varepsilon\varphi_1+\varepsilon^2\varphi_2)=\varepsilon^{-2}Q\varphi+
\varepsilon^{-1}[Q\varphi_1+Q_1\varphi]+[Q\varphi_2+Q_1\varphi_1]+\varepsilon
Q_1\varphi_2.$$

And as soon as it is equal to the right side, we obtain:
$$\left\{\begin{array}{c}
                                                Q\varphi=0, \\
                                                Q\varphi_1+Q_1\varphi=0, \\
                                                Q\varphi_2+Q_1\varphi_1=\psi, \\
                                                Q_1\varphi_2=\theta^{\varepsilon}.
                                              \end{array}
\right.$$

From the last equation we may see that the function $\varphi_2$
should be smooth enough to provide boundness of $Q_1\varphi_2.$
Moreover, from the first equation we see that the function $\varphi$
may by any function from the null-space of $Q$ and does not depend
on the variable that corresponds to the switching process.

An important condition of solvability of this problem is balance
condition
$$\Pi Q_1=0.$$ This condition means that the function $Q_1\varphi$ belongs to the range of $Q$,
thus we may solve the second equation using the potential operator,
that in inverse to $Q$ at its range
\begin{equation} \label{SP1}\varphi_1=-R_0Q_1\varphi.\end{equation}

Thus the main problem is to solve the equation
$$Q\varphi_2=\psi-Q_1\varphi_1=\psi+Q_1R_0Q_1\varphi.$$

The solvability condition for $Q$ has the view:
$$\Pi Q\Pi\varphi_2=0=\Pi\psi+\Pi Q_1R_0Q_1\Pi\varphi,$$ and we finally obtain \begin{equation} \label{SP2}\Pi\psi=-\Pi
Q_1R_0Q_1\Pi\varphi.\end{equation}

For the function $\varphi_2$ we obviously have: \begin{equation}
\label{SP3}\varphi_2=R_0[-\Pi
Q_1R_0Q_1\Pi+Q_1R_0Q_1]\varphi.\end{equation}

Equations (\ref{SP1})-(\ref{SP3}) give the solution of singular
perturbation problem.

In case of MSRE balance condition has the following view
\begin{equation}
\label{UB}\Pi S(\theta)\mathbf{1}(x)=0,\end{equation} where
$\mathbf{1}(x)=(x_1,...,x_n).$ Really, every term under the sign of
integral contains either $\int_0^{\pi}\sin^n\theta\cos \theta
d\theta=0$ or $\int_0^{2\pi}\sin\theta d\theta=0$.

\section{Main result for MSRE}

{\bf THEOREM 1.} MSRE $\xi^{\varepsilon}_t,$ converges weakly to the
Wiener process $w(t):=\xi^{0}_t$ when $\varepsilon \to 0$:

$$\xi^{\varepsilon}_t\Rightarrow \xi^{0}_t,$$
where $\xi^{0}_t\in \mathbf{R}^n$ is defined by a generator
\begin{equation} \label{Gen0}
L^{0}\varphi(x_1,...,x_n)=\frac{c^2}{n}\Delta\varphi(x_1,...,x_n),
\end{equation}
where $\Delta\varphi:=\left(\frac{\partial^2}{\partial
x_1^2}+...+\frac{\partial^2}{\partial x_n^2}\right)\varphi$.

\textbf{Remark 1:} The generator (\ref{Gen0}) corresponds to the
results, obtained in \cite{Kol2} for the spaces $\mathbf{R}^2$ and
$\mathbf{R}^3$.

To prove the Theorem we need the following Lemma.

{\bf Lemma 1.} At the perturbed test-functions
\begin{equation}
\label{phi}\varphi^{\varepsilon}(x_1,...,x_n;\theta)=\varphi(x_1,...,x_n)+\varepsilon\varphi_1(x_1,...,x_n;\theta)+\varepsilon^2\varphi_2(x_1,...,x_n;\theta),\end{equation}
that have bounded derivatives of any degree and compact support, the
operator $L^{\varepsilon}$ has asymptotic representation
$$L^{\varepsilon}\varphi^{\varepsilon}(x_1,...,x_n;\theta)=L^{0}\varphi(x_1,...,x_n)+R^{\varepsilon}(\theta)\varphi(x),$$
$$|R^{\varepsilon}(\theta)\varphi(x)|\to 0, \varepsilon\to 0,
\varphi(x)\in C^{\infty}_0(\mathbf{R}^d),$$ where $L^{0}$ is defined
in (\ref{Gen0}), $\varphi_1(x_1,...,x_n;\theta),
\varphi_2(x_1,...,x_n;\theta)$ and
$R^{\varepsilon}(\theta)\varphi(x)$ are the following:
\begin{equation} \label{L0}
L^0\Pi=-c^2\Pi S(\theta)R_0S(\theta)\Pi, \end{equation}
$$\begin{array}{c}
                                    \varphi_1=-cR_0S(\theta)\varphi, \\
                                    \varphi_2=c^2R_0S(\theta)R_0S(\theta)\varphi, \\
                                    R^{\varepsilon}(\theta)\varphi=\varepsilon
                                    c^3
                                    S(\theta)R_0S(\theta)R_0S(\theta)\varphi.
                                  \end{array}
$$

{\bf Proof.} Let's solve singular perturbation problem for the
operator (\ref{GenEp}). To do this, we study this operator at the
test-function (\ref{phi}). We have:
$$L^{\varepsilon}\varphi^{\varepsilon}(x_1,...,x_n;\theta)=[\varepsilon^{-2}Q+\varepsilon^{-1}cS(\theta)][\varphi+\varepsilon\varphi_1+\varepsilon^2\varphi_2]=
\varepsilon^{-2}Q\varphi+$$
$$\varepsilon^{-1}[Q\varphi_1+cS(\theta)\varphi]+[Q\varphi_2+cS(\theta)\varphi_1]+\varepsilon cS(\theta)\varphi_2.$$

Thus, we obtain the following equations
\begin{equation} \label{Eqn}\left\{\begin{array}{c}
                                    Q\varphi=0, \\
                                    Q\varphi_1+cS(\theta)\varphi=0, \\
                                    L^0\varphi=Q\varphi_2+cS(\theta)\varphi_1, \\
                                    R^{\varepsilon}\varphi(\theta)=\varepsilon
                                    c
                                    S(\theta)\varphi_2.
                                  \end{array}
\right.\end{equation}

From the first equation we see that $\varphi(x_1,...,x_n)$ belongs
to the null-space of $Q$. From the balance condition (\ref{UB}) easy
to see that $S(\theta)\varphi$ belongs to the range of $Q$, thus
from the second equation of the system (\ref{Eqn}) we obtain
$$\varphi_1=-cR_0S(\theta)\varphi.$$

By substitution into the third equation and using the solvability
condition we have:
$$L^0\Pi\varphi+c^2\Pi S(\theta)R_0S(\theta)\Pi\varphi=0.$$

From the last equation of (\ref{UB}):
$$R^{\varepsilon}(\theta)\varphi(x)=\varepsilon c^3S(\theta)R_0S(\theta)R_0S(\theta)\varphi(x)\rightarrow 0 \mbox{
when  } \varepsilon\rightarrow 0, \varphi(x)\in
C^{\infty}_0(\mathbf{R}^n).$$

Let's find the generator of the limit process $L^0$ by the formula
(\ref{L0}):
$$L^0\varphi=c^2\Pi S(\theta)(I-\Pi)S(\theta)\Pi\varphi=c^2\Pi S^2(\theta)\Pi\varphi-c^2\Pi S(\theta)\Pi S(\theta)\Pi\varphi.$$

The last term equals to 0 by the balance condition (\ref{UB}). Thus,
finally:
\begin{equation} \label{For}L^0=c^2\Pi S^2(\theta),\end{equation}
or
$$L^0=\frac{c^2}{N}\int_{S_n}S^2(\theta)\mu(d\theta).$$

Let's calculate the integral:
$$\frac{c^2}{N}\int_0^{\pi}...\int_0^{\pi}\int_0^{2\pi}\left[\cos^2\theta_1\frac{\partial^2}{\partial
x_1^2}+\sin^2\theta_1\cos^2\theta_2\frac{\partial^2}{\partial
x_2^2}+...+\right.$$ $$\left.
\sin^2\theta_1...\sin^2\theta_{n-2}\cos^2\theta_{n-1}\frac{\partial^2}{\partial
x_{n-1}^2}+\sin^2\theta_1...\sin^2\theta_{n-2}\sin^2\theta_{n-1}\frac{\partial^2}{\partial
x_{n}^2}+\right.$$
$$\left.\left\{\sin\theta_1\cos\theta_1\cos\theta_2\frac{\partial^2}{\partial
x_1\partial x_2}+...+\sin^2\theta_1...
\sin^2\theta_{n-2}\sin\theta_{n-1}\times\right.\right.$$
$$\left.\left.\cos\theta_{n-1}\frac{\partial^2}{\partial x_{n-1}\partial
x_n}\right\}\right]\sin^{n-2}\theta_1\sin^{n-3}\theta_2...\sin\theta_{n-2}d\theta_1...d\theta_{n-1}=$$
$$\frac{c^2}{N}\int_0^{\pi}...\int_0^{\pi}\int_0^{2\pi}\left[\cos^2\theta_1\sin^{n-2}\theta_1\sin^{n-3}\theta_2...\sin\theta_{n-2}\frac{\partial^2}{\partial
x_1^2}+\right.$$
$$\sin^{n}\theta_1\cos^2\theta_2\sin^{n-3}\theta_2...\sin\theta_{n-2}\frac{\partial^2}{\partial
x_2^2}+...+$$ $$
\sin^{n}\theta_1\sin^{n-1}\theta_2...\sin^3\theta_{n-2}\cos^2\theta_{n-1}\frac{\partial^2}{\partial
x_{n-1}^2}+$$
$$\sin^{n}\theta_1\sin^{n-1}\theta_2...\sin^2\theta_{n-2}\sin^2\theta_{n-1}\frac{\partial^2}{\partial
x_{n}^2}+$$
$$\left\{\sin^{n-1}\theta_1\cos\theta_1\cos\theta_2\sin^{n-3}\theta_2...\sin\theta_{n-2}\frac{\partial^2}{\partial
x_1\partial x_2}+...+\right.$$
$$\left.\left.\sin^{n}\theta_1\sin^{n-1}\theta_2...\sin^3\theta_{n-2}\sin\theta_{n-1}\cos\theta_{n-1}\frac{\partial^2}{\partial x_{n-1}\partial
x_n}\right\}\right]d\theta_1...d\theta_{n-1}=$$
$$\frac{c^2}{N}\int_0^{\pi}...\int_0^{\pi}\int_0^{2\pi}\left[(\sin^{n-2}\theta_1\sin^{n-3}\theta_2...\sin\theta_{n-2}-
\sin^{n}\theta_1\sin^{n-3}\theta_2...\sin\theta_{n-2})\frac{\partial^2}{\partial
x_1^2}+\right.$$
$$(\sin^{n}\theta_1\sin^{n-3}\theta_2...\sin\theta_{n-2}-\sin^{n}\theta_1\sin^{n-1}\theta_2\sin^{n-4}\theta_3...\sin\theta_{n-2})\frac{\partial^2}{\partial
x_2^2}+...+$$ $$
(\sin^{n}\theta_1\sin^{n-1}\theta_2...\sin^3\theta_{n-2}-\sin^{n}\theta_1\sin^{n-1}\theta_2...\sin^3\theta_{n-2}\sin^2\theta_{n-1})\frac{\partial^2}{\partial
x_{n-1}^2}+$$
$$\sin^{n}\theta_1\sin^{n-1}\theta_2...\sin^3\theta_{n-2}\sin^2\theta_{n-1}\frac{\partial^2}{\partial
x_{n}^2}+$$
$$\left\{\sin^{n-1}\theta_1\cos\theta_1\cos\theta_2\sin^{n-3}\theta_2...\sin\theta_{n-2}\frac{\partial^2}{\partial
x_1\partial x_2}+...+\right.$$
$$\left.\left.\sin^{n}\theta_1\sin^{n-1}\theta_2...\sin^3\theta_{n-2}\sin\theta_{n-1}\cos\theta_{n-1}\frac{\partial^2}{\partial x_{n-1}\partial
x_n}\right\}\right]d\theta_1...d\theta_{n-1}.$$

Every term in braces has a multiplier of the type
$\int_0^{\pi}\sin^{n}\theta\cos\theta d\theta=0$ or
$\int_0^{2\pi}\sin\theta\cos\theta d\theta=0$, thus the
corresponding integral equals to 0.

Integration of every correlation in parentheses gives
$$\int_0^{\pi}...\int_0^{\pi}\int_0^{2\pi}(\sin^{n}\theta_1\sin^{n-1}\theta_2...\sin^{n-k+2}\theta_{k-1}\sin^{n-k-1}\theta_k\sin^{n-k-2}\theta_{k+1}...\sin\theta_{n-2}-$$
$$\sin^{n}\theta_1\sin^{n-1}\theta_2...\sin^{n-k+2}\theta_{k-1}\sin^{n-k+1}\theta_k\sin^{n-k-2}\theta_{k+1}...\sin\theta_{n-2})d\theta_1...d\theta_{n-1}=$$
$$N\left(\frac{(2m-1)(2m-2)(2m-3)\cdot...\cdot(2m-k+1)}{2m(2m-1)(2m-2)\cdot...\cdot(2m-k+2)}-\right.$$ $$\left.\frac{(2m-1)(2m-2)(2m-3)\cdot...\cdot(2m-k)}{2m(2m-1)(2m-2)\cdot...\cdot(2m-k+1)}\right)=$$ $$N\left(\frac{2m-k+1}{2m}-\frac{2m-k}{2m}\right),
\mbox{if}\hskip 2mm n=2m$$
$$\mbox{or}   \hskip 2mm        N\left(\frac{(2m)(2m-1)(2m-2)\cdot...\cdot(2m-k+2)}{(2m+1)2m(2m-1)\cdot...\cdot(2m-k+3)}-\right.$$ $$\left.\frac{2m(2m-1)(2m-2)\cdot...\cdot(2m-k+1)}{(2m+1)2m(2m-1)\cdot...\cdot(2m-k+2)}\right)=$$ $$N\left(\frac{2m-k+2}{2m+1}-\frac{2m-k+1}{2m+1}\right),
\mbox{if} \hskip 2mm n=2m+1$$
$$=\frac{N}{n}. $$

Finally, we have:
$$L^0\varphi(x_1,...,x_n)=\frac{c^2}{n}\Delta\varphi(x_1,...,x_n).$$

Lemma is proved.

{\bf Proof of Theorem 1.} In 1 we proved that
$L^{\varepsilon}\varphi^{\varepsilon}\Rightarrow L^0\varphi$ at the
class of functions $C^{\infty}_0(\mathbf{R}^n\times\Theta)$ when
$\varepsilon\to 0$. To prove the weak convergence we need to show
the relative compactness of the family $(\xi^{\varepsilon}_t,
\theta^{\varepsilon}_t)$ in
$\mathbf{D}_{\mathbf{R}^n\times\Theta}[0,\infty)$. To do this we use
the methods proposed in \cite{EtKu, korlim, stroock}. Let's
formulate Corollary 6.1 from \cite{korlim} (see also Theorem 6.4 in
\cite{korlim}) as a Lemma.

{\bf Lemma 2.} Let the generators $L^{\varepsilon}, \varepsilon>0$
satisfy the inequalities
$$|L^{\varepsilon}\varphi(u)|<C_{\varphi}$$ for any real-valued non-negative function $\varphi\in
C^{\infty}_0(\mathbf{R}^n\times \Theta)$, where the constant
$C_{\varphi}$ depends on the norm of $\varphi$, and for
$\varphi_0(u)=\sqrt{1+u^2},$
$$L^{\varepsilon}\varphi_0(u)\leq C_l\varphi_0(u), |u|\leq l,$$ where the constant $C_l$
depends on the function $\varphi_0$, but do not depend on
$\varepsilon>0.$

Then the family $(\xi_t^{\varepsilon},\theta_t^{\varepsilon}),
t\geq0,\varepsilon>0$ is relatively compact in
$\mathbf{D}_{\mathbf{R}^n\times\Theta}[0,\infty)$.

Let's study (\ref{GenEp}) at the test-function
$\varphi^{\varepsilon}(x_1,...,x_n;\theta)=\varphi(x_1,...,x_n)+\varepsilon\varphi_1(x_1,...,x_n;\theta),$
where $\varphi_1(x_1,...,x_n;\theta)=-cR_0S(\theta)\varphi
(x_1,...,x_n).$

We have:
$$L^{\varepsilon}\varphi^{\varepsilon}(x_1,...,x_n;\theta)=\varepsilon^{-2}Q\varphi(x_1,...,x_n)+\varepsilon^{-1}[Q\varphi_1+cS(\theta)\varphi(x_1,...,x_n)]+$$ $$cS(\theta)\varphi_1(x_1,...,x_n;\theta).$$

It follows from (\ref{Eqn}) that two first terms equal to 0. Let's
estimate the last term:
$$cS(\theta)\varphi_1(x_1,...,x_n;\theta)=c^2S(\theta)R_0S(\theta)\varphi(x_1,...,x_n)=c^2S^2(\theta)\varphi(x_1,...,x_n)=$$ $$c^2\left[\cos^2\theta_1\frac{\partial^2}{\partial
x_1^2}+\sin^2\theta_1\cos^2\theta_2\frac{\partial^2}{\partial
x_2^2}+...+\right.$$ $$\left.
\sin^2\theta_1...\sin^2\theta_{n-2}\cos^2\theta_{n-1}\frac{\partial^2}{\partial
x_{n-1}^2}+\sin^2\theta_1...\sin^2\theta_{n-2}\sin^2\theta_{n-1}\frac{\partial^2}{\partial
x_{n}^2}+\right.$$
$$\left.\left\{\sin\theta_1\cos\theta_1\cos\theta_2\frac{\partial^2}{\partial
x_1\partial x_2}+...+\sin^2\theta_1...
\sin^2\theta_{n-2}\sin\theta_{n-1}\times\right.\right.$$
$$\left.\left.\cos\theta_{n-1}\frac{\partial^2}{\partial x_{n-1}\partial
x_n}\right\}\right]\varphi(x_1,...,x_n)\leq C_{1,\varphi},$$ as soon
as all the constants, functions and their derivatives are bounded.

We also have from (\ref{Eqn}):
$$L^{\varepsilon}\varphi^{\varepsilon}=L^{\varepsilon}\varphi+\varepsilon L^{\varepsilon}\varphi_1=L^{\varepsilon}\varphi+\varepsilon cL^{\varepsilon}R_0S(\theta)\varphi.$$

Thus,
$$L^{\varepsilon}\varphi=L^{\varepsilon}\varphi^{\varepsilon}-\varepsilon cL^{\varepsilon}R_0S(\theta)\varphi\leq C_{1,\varphi}-\varepsilon
C_{2,\varphi}<C_{\varphi}$$ for small $\varepsilon.$

To prove the second condition, it's enough to use the properties of
the function $\varphi_0(u)=\sqrt{1+u^2},$ namely:
$$|\varphi'_0(u)|\leq1\leq\varphi_0(u), |\varphi''_0(u)|\leq\varphi_0(u).$$

So, the proof of the second condition is similar to the previous
reasoning.

Thus, the family $(\xi_t^{\varepsilon},\theta_t^{\varepsilon})$ is
relatively compact in
$\mathbf{D}_{\mathbf{R}^n\times\Theta}[0,\infty)$.

Now we may use the following theorem (Theorem 6.6 from
\cite{korlim}).

{\bf THEOREM 2.} Let random evolution with Markov switching
$(\xi^{\varepsilon}(t),x^{\varepsilon}(t))\in
\mathbf{D}_{\mathbf{R}^n\times E}[0,\infty)$ satisfies the
conditions:

\textbf{C1:} The family of processes $(\xi^{\varepsilon}(t), t\geq
0, \varepsilon
> 0$ is relatively compact.

\textbf{C2:} There exists a family of test-functions
$\varphi^{\varepsilon}(u,x)\in C^3_0(\mathbf{R}^n\times E)$ such
that
$$\lim\limits_{\varepsilon\to
0}\varphi^{\varepsilon}(u,x)=\varphi(u)$$ uniformly by $u,x.$

\textbf{C3:} The following uniform convergence is true
$$\lim\limits_{\varepsilon\to
0}L^{\varepsilon}\varphi^{\varepsilon}(u,x)=L\varphi(u)$$ uniformly
by $u,x.$

The family $L^{\varepsilon}\varphi^{\varepsilon}, \varepsilon>0$ is
uniformly bounded, moreover $L^{\varepsilon}\varphi^{\varepsilon}$
and $L\varphi$ belong to $C(\mathbf{R}^n\times E).$

\textbf{C4:} Convergence by probability of initial values
$$\xi^{\varepsilon}(0)\rightarrow \widehat{\xi}(0),
\varepsilon\rightarrow0,$$ and
$$\sup\limits_{\varepsilon>0}\mathbf{E}|\xi^{\varepsilon}(0)|\leq
C<+\infty$$ is true.

Then we have the weak convergence
$$\xi^{\varepsilon}(t)\Rightarrow\widehat{\xi}(t), \varepsilon\rightarrow 0.$$

According to Theorem 2, we may confirm the weak convergence in
$\mathbf{D}_{\mathbf{R}^n}[0,\infty)$
$$\xi^{\varepsilon}_t\Rightarrow\xi^{0}_t.$$
Really, all the conditions are satisfied. Namely, the family of
processes is relatively compact, the generators at the
test-functions from the class
$C^{\infty}_0(\mathbf{R}^n\times\Theta)$ converge, initial
conditions for the limit and pre-limit processes are equal.

Theorem is proved.

\section{Description of MNRE}

In $\mathbf{R}^n$ we study the same particle as in the Section 2,
but its velocity is equal to
$v(\theta)=c(\theta)\varepsilon^{-1}+c_1(\theta)$, where
$\varepsilon\to 0 $ $ (\varepsilon>0)$ is the small parameter,
functions $c(\theta), c_1(\theta)$ are bounded.

\textbf{Definition 2:} \textit{Markov non-symmetry random evolution}
(MNRE) is the process $\widetilde{\xi}_t^{\varepsilon}\in
\mathbf{R}^n$, of the following view:
$$\widetilde{\xi}_t^{\varepsilon}:=x+\int_0^tv(\theta_{\tau}^{\varepsilon})s(\theta_{\tau}^{\varepsilon})d\tau.$$

Our aim is to prove the weak convergence of MNRE to a diffusion
process with drift when $\varepsilon\to 0$.

Two-component Markov process $(\widetilde{\xi}_t^{\varepsilon},
\theta_t^{\varepsilon})$ at the test-functions
$\varphi(x_1,...,x_n;\theta)\in C^{\infty}_0(\mathbf{R}^n\times
\Theta)$ is defined by a generator (see, e.g. \cite{Pin})
\begin{equation} \label{GenEp2} L^{\varepsilon}\varphi(x_1,...,x_n;\theta):=\lambda Q\varphi(\cdot;\theta)+v(\theta) S(\theta)\varphi(x_1,...,x_n;\cdot)=$$ $$\varepsilon^{-2}Q\varphi(\cdot;\theta)+\varepsilon^{-1}c(\theta)S(\theta)\varphi(x_1,...,x_n;\cdot)+c_1(\theta)S(\theta)\varphi(x_1,...,x_n;\cdot),\end{equation}
where
$$S(\theta)\varphi(x_1,...,x_n;\cdot):=-\left(s(\theta),\nabla\right)\varphi(x_1,...,x_n;\cdot).$$

An important condition that allows to confirm weak convergence is
balance condition
\begin{equation}
\label{UB2}\Pi c(\theta)S(\theta)=0.\end{equation}

\textbf{Remark 2:} This is the last condition that defines symmetry
or non-symmetry of the process. In case of MSRE the balance
condition (\ref{UB2}) is true, and $c_1(\theta)\equiv 0.$ In case of
MNRE non-symmetry of the process is caused by the condition:
\begin{equation} \label{NUB}\Pi c_1(\theta)S(\theta)=(d,\nabla)\neq
0.\end{equation}

\textbf{Example 1:} Condition (\ref{UB2}) may be satisfied for
different functions $c(\theta)$. Namely, in case of MSRE
$c(\theta)=c=const$. Then every term under the integral contains
either $\int_0^{\pi}\sin^n\theta\cos \theta d\theta=0$ or
$\int_0^{2\pi}\sin\theta d\theta=0$.

Another variant of the function $c(\theta)$, is $c(\theta)=sin
\theta_1$. Really, under the integral we obtain terms, analogical to
the previous ones. We note that the dimension of the space i this
case should be more then 2, because in $\mathbf{R}^2$ this function
does not preserve the symmetry.

\textbf{Example 2:} Non-symmetry condition (\ref{NUB}) is also
satisfied for different functions $c_1(\theta)$. For example, in
$\mathbf{R}^2$ for $c_1(\theta)=sin \theta$ we obtain:
$$\frac{1}{2\pi}\int_0^{2\pi}sin\theta\left[cos\theta\frac{\partial}{\partial x_1}+sin\theta\frac{\partial}{\partial x_2}\right]d\theta
=\frac{1}{2}\frac{\partial}{\partial x_2}.$$

Another example is the function in $\mathbf{R}^n$:
$$c_1(\theta)=\left\{\begin{array}{c}
                                       c_1, \theta_{n-1}\in [\pi,2\pi), \\
                                       0, \theta_{n-1}\in [0,\pi).
                                     \end{array}
\right.$$ Again, all terms under the integral, except the last one
contain $\int_0^{\pi}\sin^n\theta\cos \theta d\theta=0$, so only one
term is not trivial
$$\frac{1}{N}c_1\int_0^{\pi}...\int_0^{\pi}\int_{\pi}^{2\pi}sin^{n-1}\theta_1sin^{n-2}\theta_2...sin^{2}\theta_{n-1}sin\theta_{n-1}d\theta_1...d\theta_{n-1}\frac{\partial}{\partial x_n}.$$
By simple calculations, we obtain:
$$\Pi c_1(\theta)S(\theta)=\left\{\begin{array}{c}
                                              -\frac{c_1}{2}\frac{3\cdot 5\cdot...\cdot(n-2)}{2\cdot4\cdot...\cdot(n-1)}\frac{\partial}{\partial x_n}, n=2m+1, \\
                                              -\frac{c_1}{\pi}[1\cdot2\cdot...(n-2)]\frac{\partial}{\partial x_n},
                                              n=2m.
                                            \end{array}
\right.$$

Thus we have a wide range of functions that preserve or, on the
contrary, do not preserve symmetry. So, we may define the velocity
of random evolution in different ways.

\section{Main result for MNRE}

{\bf THEOREM 3.} MNRE $\widetilde{\xi}^{\varepsilon}_t,$ converges
weakly to the process $\widetilde{\xi}^{0}_t$ when $\varepsilon \to
0$:

$$\widetilde{\xi}^{\varepsilon}_t\Rightarrow \widetilde{\xi}^{0}_t.$$

The limit process $\widetilde{\xi}^{0}_t\in \mathbf{R}^n$ is defined
by a generator
\begin{equation} \label{Gen20}
L^{0}\varphi(x_1,...,x_n)=(d,\nabla)\varphi(x_1,...,x_n)+(\sigma^2,\Delta)\varphi(x_1,...,x_n),
\end{equation}
where $\Delta\varphi(x_1,...,x_n):=\left(\frac{\partial^2}{\partial
x_1^2}+...+\frac{\partial^2}{\partial
x_n^2}\right)\varphi(x_1,...,x_n),$
$(d,\nabla):=-\frac{1}{N}\times\int_{S_n}c_1(\theta)(s(\theta),\nabla)\mu(d\theta),$
$(\sigma^2,\Delta):=\frac{1}{N}\int_{S_n}c^2(\theta)(s(\theta),\nabla)^2\mu(d\theta).$

To prove the Theorem, we need the following Lemma.

{\bf Lemma 3.} At the perturbed test-functions
\begin{equation}
\label{phi2}\varphi^{\varepsilon}(x_1,...,x_n;\theta)=\varphi(x_1,...,x_n)+\varepsilon\varphi_1(x_1,...,x_n;\theta)+\varepsilon^2\varphi_2(x_1,...,x_n;\theta),\end{equation}
that have bounded derivatives of any degree and compact support, the
operator $L^{\varepsilon}$ has asymptotic representation
$$L^{\varepsilon}\varphi^{\varepsilon}(x_1,...,x_n;\theta)=L^{0}\varphi(x_1,...,x_n)+R^{\varepsilon}(\theta)\varphi(x),$$ $$|R^{\varepsilon}(\theta)\varphi(x)|\to 0, \varepsilon\to 0, \varphi(x)\in C^{\infty}(\mathbf{R}^d),$$
where $L^{0}$ is defined in (\ref{Gen20}),
$\varphi_1(x_1,...,x_n;\theta), \varphi_2(x_1,...,x_n;\theta)$ and
$R^{\varepsilon}(\theta)\varphi(x)$ are defined by equalities:
\begin{equation} \label{L20}
L^0\Pi=-\Pi c(\theta)S(\theta)R_0c(\theta)S(\theta)\Pi+\Pi
c_1(\theta)S(\theta)\Pi\varphi,
\end{equation}
$$\begin{array}{c}
                                    \varphi_1=-R_0c(\theta)S(\theta)\varphi, \\
                                    \varphi_2=R_0c(\theta)S(\theta)R_0c(\theta)S(\theta)\varphi, \\
                                    R^{\varepsilon}(\theta)\varphi=\left\{\varepsilon [c(\theta)S(\theta)R_0c(\theta)S(\theta)R_0c(\theta)S(\theta)+c_1(\theta)S(\theta)R_0c(\theta)S(\theta)]+\right. \\ \left.\varepsilon^2c_1(\theta)S(\theta)R_0c(\theta)S(\theta)R_0c(\theta)S(\theta)\right\}\varphi.
                                  \end{array}
$$

{\bf Proof.} We solve singular perturbation problem for the operator
(\ref{GenEp2}). To do this, we study this operator at the
test-function (\ref{phi2}). So, we have:
$$L^{\varepsilon}\varphi^{\varepsilon}(x_1,...,x_n;\theta)=[\varepsilon^{-2}Q+\varepsilon^{-1}c(\theta)S(\theta)+c_1(\theta)S(\theta)][\varphi+\varepsilon\varphi_1+\varepsilon^2\varphi_2]=$$
$$
\varepsilon^{-2}Q\varphi+\varepsilon^{-1}[Q\varphi_1+c(\theta)S(\theta)\varphi]+[Q\varphi_2+c(\theta)S(\theta)\varphi_1+c_1(\theta)S(\theta)\varphi]+$$
$$\varepsilon [c(\theta)S(\theta)\varphi_2+c_1(\theta)S(\theta)\varphi_1]+\varepsilon^2c_1(\theta)S(\theta)\varphi_2.$$

Thus, we obtain the following equations:
\begin{equation} \label{Eqn2}\left\{\begin{array}{c}
                                    Q\varphi=0, \\
                                    Q\varphi_1+c(\theta)S(\theta)\varphi=0, \\
                                    L^0\varphi=Q\varphi_2+c(\theta)S(\theta)\varphi_1+c_1(\theta)S(\theta)\varphi, \\
                                    R^{\varepsilon}\varphi(\theta)=\varepsilon [c(\theta)S(\theta)\varphi_2+c_1(\theta)S(\theta)\varphi_1]+\varepsilon^2c_1(\theta)S(\theta)\varphi_2.
                                  \end{array}
\right.\end{equation}

According to the first one $\varphi(x_1,...,x_n)$ belongs to the
null-space of $Q$. From the balance condition (\ref{UB2}) we see
that $c(\theta)S(\theta)\varphi$ belongs to the range of $Q$, thus
from the second equation of (\ref{Eqn2}) we have
$$\varphi_1=-R_0c(\theta)S(\theta)\varphi.$$

By substitution into the third equation, and using the solvability
condition, we obtain:
$$L^0\Pi\varphi+\Pi
c(\theta)S(\theta)R_0c(\theta)S(\theta)\Pi\varphi-\Pi
c_1(\theta)S(\theta)\Pi\varphi=0.$$

From the last equation:
$$R^{\varepsilon}(\theta)\varphi(x)=\left\{\varepsilon [c(\theta)S(\theta)R_0c(\theta)S(\theta)R_0c(\theta)S(\theta)+c_1(\theta)S(\theta)R_0c(\theta)S(\theta)]+\right.$$
$$\left.\varepsilon^2c_1(\theta)S(\theta)R_0c(\theta)S(\theta)R_0c(\theta)S(\theta)\right\}\varphi(x)\rightarrow 0 \mbox{
при  } \varepsilon\rightarrow 0, \varphi(x)\in
C^{\infty}_0(\mathbf{R}^n).$$

Let's calculate the operator $L^0$ by the formula (\ref{L20}):
$$L^0\varphi=\Pi c(\theta)S(\theta)(I-\Pi)c(\theta)S(\theta)\Pi\varphi+\Pi
c_1(\theta)S(\theta)\Pi\varphi=\Pi
c^2(\theta)S^2(\theta)\Pi\varphi-$$ $$\Pi c(\theta)S(\theta)\Pi
c(\theta)S(\theta)\Pi\varphi+\Pi c_1(\theta)S(\theta)\Pi\varphi.$$

The second term equals 0 by the balance condition (\ref{UB2}), the
last one is not equal to 0 by (\ref{NUB}). Thus, we finally have:
\begin{equation} \label{For2}L^0=\Pi c^2(\theta)S^2(\theta)+\Pi c_1(\theta)S(\theta).\end{equation}

Using the view of $S(\theta)$, we may write: $$\Pi
c_1(\theta)S(\theta)=-\frac{1}{N}\int_{S_n}c_1(\theta)(s(\theta),\nabla)\mu(d\theta)=:(d,\nabla),$$
$$\Pi c^2(\theta)S^2(\theta)=\frac{1}{N}\int_{S_n}c^2(\theta)(s(\theta),\nabla)^2\mu(d\theta)=:(\sigma^2,\Delta).$$

Lemma is proved.

{\bf Proof of theorem 3.} We proved in Lemma 3 that
$L^{\varepsilon}\varphi^{\varepsilon}\Rightarrow L^0\varphi$ at the
class $C^{\infty}_0(\mathbf{R}^n\times\Theta)$ when $\varepsilon\to
0$. To prove the weak convergence we should show relative
compactness of the family $(\widetilde{\xi}^{\varepsilon}_t,
\theta^{\varepsilon}_t)$ in
$\mathbf{D}_{\mathbf{R}^n\times\Theta}[0,\infty)$. To do this, we
use Lamma 2.

Let's study the operator (\ref{GenEp2}) at the test-function
$\varphi^{\varepsilon}(x_1,...,x_n;\theta)=\varphi(x_1,...,x_n)+\varepsilon\varphi_1(x_1,...,x_n;\theta),$
where $\varphi_1(x_1,...,x_n;\theta)=-R_0c(\theta)S(\theta)\varphi
(x_1,...,x_n).$

We have:
$$L^{\varepsilon}\varphi^{\varepsilon}(x_1,...,x_n;\theta)=\varepsilon^{-2}Q\varphi(x_1,...,x_n)+\varepsilon^{-1}[Q\varphi_1+c(\theta)S(\theta)\varphi(x_1,...,x_n)]+$$ $$[c(\theta)S(\theta)\varphi_1(x_1,...,x_n;\theta)+c_1(\theta)S(\theta)\varphi(x_1,...,x_n;\theta)]+\varepsilon c_1(\theta)S(\theta)\varphi_1(x_1,...,x_n;\theta).$$

It follows from (\ref{Eqn2}) that the first two terms equals to 0.
Let' estimate the third term:
$$c(\theta)S(\theta)\varphi_1(x_1,...,x_n;\theta)+c_1(\theta)S(\theta)\varphi(x_1,...,x_n;\theta)=$$ $$c(\theta)S(\theta)R_0c(\theta)S(\theta)\varphi(x_1,...,x_n)+c_1(\theta)S(\theta)\varphi(x_1,...,x_n;\theta)=$$ $$c^2(\theta)S^2(\theta)\varphi(x_1,...,x_n)+c_1(\theta)S(\theta)\varphi(x_1,...,x_n;\theta)=$$ $$c^2(\theta)\left[\cos^2\theta_1\frac{\partial^2}{\partial
x_1^2}+\sin^2\theta_1\cos^2\theta_2\frac{\partial^2}{\partial
x_2^2}+...+\right.$$ $$\left.
\sin^2\theta_1...\sin^2\theta_{n-2}\cos^2\theta_{n-1}\frac{\partial^2}{\partial
x_{n-1}^2}+\sin^2\theta_1...\sin^2\theta_{n-2}\sin^2\theta_{n-1}\frac{\partial^2}{\partial
x_{n}^2}+\right.$$
$$\left.\left\{\sin\theta_1\cos\theta_1\cos\theta_2\frac{\partial^2}{\partial
x_1\partial x_2}+...+\sin^2\theta_1...
\sin^2\theta_{n-2}\sin\theta_{n-1}\times\right.\right.$$
$$\left.\left.\cos\theta_{n-1}\frac{\partial^2}{\partial x_{n-1}\partial
x_n}\right\}\right]\varphi(x_1,...,x_n)+c_1(\theta)\left[\cos\theta_1\frac{\partial}{\partial
x_1}+\sin\theta_1\cos\theta_2\frac{\partial}{\partial
x_2}+...+\right.$$ $$\left.
\sin\theta_1...\sin\theta_{n-2}\cos\theta_{n-1}\frac{\partial}{\partial
x_{n-1}}+\sin\theta_1...\sin\theta_{n-2}\sin\theta_{n-1}\frac{\partial}{\partial
x_{n}}\right]\varphi(x_1,...,x_n)\leq C_{1,\varphi},$$ as soon as
all the constants, functions and their derivatives are bounded.

The last term may be estimated analogically.

From (\ref{Eqn2}) we also have:
$$L^{\varepsilon}\varphi^{\varepsilon}=L^{\varepsilon}\varphi+\varepsilon L^{\varepsilon}\varphi_1=L^{\varepsilon}\varphi+\varepsilon L^{\varepsilon}R_0с(\theta)S(\theta)\varphi.$$

Thus,
$$L^{\varepsilon}\varphi=L^{\varepsilon}\varphi^{\varepsilon}-\varepsilon cL^{\varepsilon}R_0S(\theta)\varphi\leq C_{1,\varphi}-\varepsilon
C_{2,\varphi}<C_{\varphi}$$ for $\varepsilon$ that are small enough.

To prove the second condition of Lemma 2, we use the following
properties of the function $\varphi_0(u)=\sqrt{1+u^2}$:
$$|\varphi'_0(u)|\leq1\leq\varphi_0(u), |\varphi''_0(u)|\leq\varphi_0(u).$$

So, the proof of the second condition is similar to the previous
reasoning.

Thus, the family of the processes
$(\widetilde{\xi}_t^{\varepsilon},\theta_t^{\varepsilon})$ is
relatively compact in
$\mathbf{D}_{\mathbf{R}^n\times\Theta}[0,\infty)$.

Using Theorem 2 we confirm the weak convergence
 in $\mathbf{D}_{\mathbf{R}^n}[0,\infty)$:
$$\widetilde{\xi}^{\varepsilon}_t\Rightarrow\widetilde{\xi}^{0}_t.$$
Really, all the conditions are satisfied. Namely, the family of
processes is relatively compact, the generators at the
test-functions from the class
$C^{\infty}_0(\mathbf{R}^n\times\Theta)$ converge, initial
conditions for the limit and pre-limit processes are equal.

Theorem is proved.

\textbf{Example 3:} We study one more variant of evolution in
$\mathbf{R}^2$.

Let $$c(\theta)=\left\{\begin{array}{c}
                                       1, \theta=0, \\
                                       1, \theta=\pi,
                                     \end{array}
\right.$$ and $$c_1(\theta)=1, \theta=\frac{\pi}{2}.$$ In other
cases both functions equal to 0.

The balance condition (\ref{UB2}) is true for $c(\theta)$, on the
contrary, condition (\ref{NUB}) is true for $c_1(\theta)$.

The limit generator (\ref{Gen20}) has the view
$$L^{0}\varphi(x_1,x_2)=\frac{1}{2\pi}\frac{\partial}{\partial x_2}\varphi(x_1,x_2)+\frac{1}{\pi}\frac{\partial^2}{\partial x^2_1}\varphi(x_1,x_2).$$

Thus, the limit process has two parts - the drift with velocity
$\frac{1}{2\pi}$ in direction of $x_2$ coordinate and diffusion part
similar to the limit process in M.Kac model \cite{Kac} in
one-dimensional subspace, corresponding to $x_1$ coordinate.

\end{document}